\documentclass[10pt]{article}
\usepackage{graphicx}
\newcommand{\ds }{ \displaystyle }                                          
\newcommand{\bd }[1]{ {\mbox{\boldmath $#1$}} }                             
\newcommand{\spc}{ \,\,\,\, }
\newcommand{\R}{ {I\!\!R} }
\newcommand{\Z}{ {Z} }
\newcommand{\N}{ {I\!\!N} }
\newcommand{\fns} { \footnotesize }

\newcommand{\mP} {\mathcal{P}}
\newcommand{\mA} {\mathcal{A}}
\newcommand{\mB} {\mathcal{B}}

\newtheorem{theorem}{Theorem}
\newenvironment{remark}
  {\noindent \slshape Remark: \begin{quote} \small \itshape}{\end{quote}}

\DeclareGraphicsExtensions{.jpg, .eps}
\begin{document}

\begin{center}
Monte Carlo Random Walk Simulations Based on Distributed Order 
Differential Equations
\vspace*{0.5cm}

{\small
Erik Andries \\	
Department of Pathology \& Department of Mathematics and Statistics\\
The University of New Mexico, Albuquerque, New Mexico 87131\\
andriese@unm.edu\\~\\
Sabir Umarov \\
Department of Mathematics and Mechanics\\
The National University of Uzbekistan, Tashkent, Uzbekistan\\
sabir@math.unm.edu\\~\\
Stanly Steinberg 
\footnote{Partially supported by NIH grant P20 GMO67594.}\\                     
Department of Mathematics and Statistics\\
The University of New Mexico, Albuquerque, New Mexico 87131\\
stanly@math.unm.edu
}
\end{center}
%%%%%%%%%%%%%%%%%%%%%%%%%%%%%%%%%%%%%%%%%%%%%%%%%%%%%%%%%%%%%%%%%%%%%%%%%%%%%%%%
\begin{abstract}
In this paper the multi-dimensional random walk models governed by
distributed fractional order differential equations and multi-term fractional order
differential equations are constructed. The scaling limits of
these random walks to a diffusion process in the sense of distributions is proved.
Simulations based upon multi-term fractional order
differential equations are performed.
\end{abstract}

\vskip 6pt



\vskip 3pt

\section{Introduction}

\subsection{Motivation.}
In this paper we study simulation models based on distributed  
order differential equations, which we will call DODE simulations. 
This type of simulation reflects the rich structure of diffusion media, in 
which a several diffusion modes are possible. Diffusion processes with 
complex and changing modes are ubiquitous in nature (see, 
\cite{BouchaudGeorges90,ChechkinGonchar,MK,UchaykinZolotarev,Zaslavsky} and references 
therein). One of the motivations for conducting
DODE simulations is to model the 
movement of proteins on the cell membrane. Numerous experiments 
\cite{Edidin,GhoshWebb,Kusumi,Saxton01,SaxtonJacobson97} show that 
macromolecule movement through the cell membrane is distinct from Brownian 
motion. Saxton and Jacobson \cite{SaxtonJacobson97} noted that practically 
all experimental results show apparent transitions among modes of motion.

The governing
equation, which we take as a basis for our simulation models,
in general form, is
distributed space fractional order differential equation
\begin{equation}
  \label{dode0}
  D_{\ast}^{\beta} u(t,x) = \int_0^2 a(\alpha)
  D_0^{\alpha} u(t,x) d \alpha, \, t>0, \, x \in {\R}^N,
\end{equation}
where $0<\beta \leq 1$, $D_{\ast}^{\beta}$ is the Caputo fractional order 
derivative \cite{Caputo,GLU}, $D_0^{\alpha}=(-\Delta)^{{\alpha \over 2}}$ is the space fractional 
order (pseudo-differential) operator with the symbol $|\xi|^{\alpha}$. 
Note that $D_0^{\alpha}$ can be written in the form of hypersingular 
integral as well \cite{SKM}. The function $a(\alpha)$ is a positive 
integrable function (or positively defined distribution). Depending on 
$a(\alpha)$, (\ref{dode0}) may become a multi-term fractional order 
differential equation, which can possibly describe the existence of 
a finite number of diffusion regimes.
Although, the distributed order differential operators were first mentioned by
\cite{Caputo1,Caputo2} in the 1960s, the intensive study of models 
based on the distributed order differential 
equations has been started recently 
\cite{BagleyTorvic,Chechkin,Diethelm,LorenzoHartley,MMM,UG05-2,UmarovSteinberg}. 

The present report is organized as follows. In Section 2, we briefly recall the 
theoretic platform of the construction of the DODE simulation models announced 
in \cite{UmarovSteinberg}. In Section 3 we analyze 
the difference schemes associated with the
DODE models, and in Sections 4 and 5 we 
construct random walk models and simulations based on the transition 
probabilities introduced in the previous sections.

\subsection{Notation.}
In this paper, ${\R^N}$ is the $N$-dimensional Euclidean space with 
coordinates ${x} =(x_1,...,x_N)$ while ${\Z^N}$ is the $N$-dimensional 
integer-valued lattice with the lattice nodes being given by the multi-index 
notation $j=(j_1,...,j_N)$.  The letters $i$, $j$ and $k$ will be exclusively 
used for the multi-indexing of lattice nodes. We denote by $x_j = 
(h_{j_1},...,h_{j_N}), j \in {\Z}^N$, the nodes of the uniform $h$-lattice 
${\Z}_h^N$ which is defined as $(h {\Z})^N$ with $h$ being the distance 
between any two lattice nodes. We introduce a spatial grid $\{x_j=jh, j \in 
{\Z}^N \}$, with $ h>0$ and a temporal grid $\{t_n=n \tau, n=0,1,2,... \}$ 
with a fixed stepsize $\tau>0$.  Furthermore, let $u^n_j$ denote the 
discretization of the function 
$u(t,x)$ on the spatial and temporal grid at $x=x_j$ and $t=t_n$, i.e  
$u^n_j = u(t_n,x_j)$.

\section{Markovian random walks associated with the DODE}
\subsection{Particle jumps.}
Assume ${\bf X}$ to be a N-dimensional random vector \cite{MS} whose values 
range in ${\Z}^N$. Let a sequence of random vectors ${\bf X}_1, {\bf X}_2, 
...$ also be N-dimensional independent identically distributed random 
vectors, all having the same probability distribution. Consider the sequence 
of random vectors
%-------------------------------------------------------------------
\[
  {\bf S}_n = h{\bf X}_1+h{\bf X}_2+...+h{\bf X}_n, n=1,2,...
\]
%-------------------------------------------------------------------
taking ${\bf S}_0 ={\bf 0} =(0,\ldots,0) \in {\Z}^N_h$ for convenience. 
We interpret ${\bf X}_1, {\bf X}_2, ...$,
as a sequence of particle jumps
starting time $t=t_0=0$.
At time $t=t_n$, the particle takes a jump $h{\bf X}_n$ from 
${\bf S}_{n-1}$ to ${\bf S}_n$.
If $u^n_j=u(t_n,x_j)$ is the probability of 
a particle being at location $x_j$ at time $t_n$ and,
taking into account the
recursion ${\bf S}_{n+1}={\bf S}_n + h{\bf X}_{n+1}$, we have
%-------------------------------------------------------------------
\begin{equation}
  \label{recursion}
  u^{n+1}_j = \sum_{k \in {\Z}^N}p_k u^n_{j-k}, j \in {\Z}^N, \,\,
  n=0,1,...
\end{equation}
%-------------------------------------------------------------------
where the coefficients $p_k, \, k \in {\Z}^N$ are called the transition 
probabilities.  
The convergence of the sequence ${\bf S}_n$ when 
$n \rightarrow \infty$ 
means convergence of the discrete probability law (probability mass 
function)
$(u^n_j)_{j\in {\Z}^N}$, properly rescaled as explained
below, to the probability law with a density $u(t,x)$ in the
sense of distributions (in law). This is equivalent to the locally
uniform convergence of the corresponding characteristic functions
(see for details \cite{MS}). This idea is used in \cite{UG05-1,UmarovSteinberg}
to prove the convergence of the
sequence of characteristic functions of the
corresponding random walks to the
fundamental solution of distributed order diffusion equations.

\subsection{Markovian transition probabilities.}
Let the transition probabilities in Eq.(\ref{recursion}) take the form
%-------------------------------------------------------------------
\begin{equation}
  p_{k}  =  \tau q_k(\alpha,h), \,\, k \neq 0,
  \label{eq:p_k}
\end{equation}
%-------------------------------------------------------------------
where
%-------------------------------------------------------------------
\begin{equation}   
  q_k(\alpha,h) =  \int_0^2 
                    \left[
                    \frac{ a(\alpha) b(\alpha) }
                         { \vert k \vert^{N+\alpha} h^{\alpha} }
                    \right] d \alpha,
  \spc \mbox{and} \spc
  b(\alpha) = 
  \frac{
        \left[ 
          \Gamma \left( 1 + \frac{\alpha}{2} \right)
        \right]^2
        sin \left( \frac{\alpha}{2} \pi \right)
       }{
        \pi^2 2^{N-\alpha-1}
       }.
  \label{eq:q_k} 
\end{equation}
%-------------------------------------------------------------------
The transition probability $p_0$ can then be defined as 
%-------------------------------------------------------------------
\begin{equation}
  p_0 = 1 - \sum_{k \neq 0} p_k = 1 - \tau q_0(\alpha,h), 
  \label{eq:p_0}
\end{equation}
%-------------------------------------------------------------------
where
%-------------------------------------------------------------------
\begin{equation}
  q_0(\alpha,h)  =   
  \sum_{k \neq 0} q_k(\alpha,h) =  
  \sum_{k \neq 0} \int_0^2 \left[
  \frac{ a(\alpha) b(\alpha) }{ \vert k \vert^{N+\alpha} h^{\alpha} }
  \right] d \alpha, 
  \label{eq:q_0}
\end{equation}
%-------------------------------------------------------------------
Assuming that the condition 
$0 < \tau q_0(\alpha,h) \leq 1$ is fulfilled,   
the transition probabilities then satisfy the following properties:
%-------------------------------------------------------------------
\begin{enumerate}
  \item  ${\displaystyle \sum_{k \in {\Z}^N} p_k} = 1;$   
  \item  $p_k \geq 0,{k \in {\Z}^N}.$
\end{enumerate}
%-------------------------------------------------------------------
Note that the non-negativity condition\footnote{This condition 
is equivalent to the
stability condition of finite-difference schemes giving the
usual stability condition if $a(\alpha)=\delta(\alpha-2)$.}
in property 2 is linked with the Riemann zeta-function. 
Indeed, introduce the function
%-------------------------------------------------------------------
\begin{equation}
  {\mathcal{R}} (\alpha) = 
  \sum_{k \neq 0} \frac{1}{|k|^{N+\alpha}} =
  \sum_{m=1}^{\infty} \frac{M_m}{m^{N+\alpha}}, \spc 0<\alpha \leq 2,
\end{equation}
%-------------------------------------------------------------------
where $M_m = \sum_{|k|=m}1.$ In the one-dimensional case
${\mathcal{R}}(\alpha)= 2 \zeta (1+\alpha)$, where $\zeta (z)$ is the
Riemann zeta-function.
Then the nonnegativity condition $0 < p_0 \leq 1$ can be rewritten as 
%-------------------------------------------------------------------
\begin{equation}
  \label{eq:condfortrpr1}
  \tau q_0(\alpha, h) = \tau \int_0^2 \left[
  \frac{a(\alpha)b(\alpha){\mathcal{R}}(\alpha)}{h^{\alpha}} 
  \right] d \alpha  \leq 1.
\end{equation}
%-------------------------------------------------------------------
It follows from this condition that $h \rightarrow 0$ yields
$\tau \rightarrow 0$. This, in turn,
yields $t/\tau \rightarrow \infty$ for any finite $t.$

%%%%%%%%%THEOREM-1%%%%%%%%%%%
\begin{theorem}
\label{theorem:dode}
Let ${\bf X}$ be a random vector with the transition probabilities
$p_{k} = P({\bf X}=x_k), k \in {\Z}^N,$ defined in 
Eq.(\ref{eq:p_k}) and Eq.(\ref{eq:p_0})
which satisfy properties 1 and 2. 
Then the sequence of random vectors 
${\bf S}_{n}=h{\bf X}_{1}+ ... + h{\bf X}_{n},$ 
converges as $n \rightarrow \infty$ in law to
the random vector whose probability density function is the
fundamental solution of the distributed space fractional order
differential equation (\ref{dode0}) with $\beta=1$.
%(\ref{cc}), i.e. $G(t,x)$ defined in (\ref{fs}).
\end{theorem}
%%%%%%%%%%%%%%%%%

Note, for the simulations used in this paper, it is important to use 
the multi-term analog of this theorem.  Assuming that
%-------------------------------------------------------------------
\begin{equation}   
  \label{discret} 
  a(\alpha)=\sum_{m=1}^{M} a_m \delta (\alpha - \alpha_m), 
  \spc
  0 < \alpha_1 < \cdots < \alpha_M \leq 2,
\end{equation}
%-------------------------------------------------------------------
with positive constants $a_m$, we get a multiterm DODE
%-------------------------------------------------------------------
\begin{equation}
  \label{multiterm}
  D_{\ast}^{\beta} u(t,x) = \sum_{m=1}^M a_m
  D_0^{\alpha_m}u(t,x), \spc 
  t>0, \, x \in {\R}^N.
\end{equation}
%-------------------------------------------------------------------
Also note that the coefficients $q_k(\alpha,h)$ in Eq.(\ref{eq:q_k})
and Eq.(\ref{eq:q_0}) become multi-term as well: 
%-------------------------------------------------------------------
\[ 
  q_k(\alpha,h) =   
  {\ds \sum_{m=1}^M 
  \left[
    \frac{ a_m b(\alpha_m) }
         { \vert k \vert^{N+\alpha_m} h^{\alpha_m} }
  \right]}, \,\, k \neq 0, 
  \quad 
  q_0 = {\ds \sum_{k \neq 0} q_k}.  
\]
%-------------------------------------------------------------------

%%%%%%%%%%%%%%%( THEOREM 2
\begin{theorem}
\label{theorem:multiterm}
Let the transition probabilities
$p_{k} = P({\bf X}=x_k), k \in {\Z}^N,$ of the
random vector ${\bf X}$ be given as follows:
%----------------------------------------------------------
\begin{equation}
  p_k = \tau q_k(\alpha,h) 
 \spc \mbox{and} \spc
  p_0 = 1 - \tau q_0(\alpha,h)
  \label{eq:markovian_transprob} 
\end{equation}
%----------------------------------------------------------
where $a(\alpha) = \sum_{m=1}^M a_m \delta(\alpha-\alpha_m)$. 
Assume
%----------------------------------------------------------
\[
  \tau
  \sum_{m=1}^M 
  \frac{ a_m b(\alpha_m){\mathcal{R}}(\alpha_m) } 
       { h^{\alpha_m} }
  \leq 1.
\]
%----------------------------------------------------------
Then the sequence of random vectors 
${\bf S}_{n}=h{\bf X}_{1}+ ... + h{\bf X}_{n},$ 
converges as $n \rightarrow \infty$ in law to
the random vector whose probability density function is the
fundamental solution of the multiterm fractional order differential
equation (\ref{multiterm}) with $\beta=1$.
\end{theorem}
%%%%%%%%%%%%%%%( END OF THEOREM 2

\begin{remark}
As we noted above these results were announced in \cite{UmarovSteinberg}. The more general case of 
these theorems corresponding to a fractional $\beta \in (0,1)$ can be 
obtained introducing a positive waiting time distribution and corresponding 
iid random variables \cite{GM,MMM}. We do not describe this case in this 
paper. We note only that the general case is studied by applying 
a general finite-difference approach and that this 
general difference scheme is stable under some 
condition and has a unique solution.
\end{remark}    
\vspace{0.5cm}

\section{Generalized Transition Probabilities for the DODE}

The set of grid points in $\Z^N_h$
used to update $u$ at time $t=t_{n+1}=(n+1) \tau$
is called the stencil.
In this section, we start from stating the values of the transition
probabilities associated with the stencil for the
discretization of the particular space-time-fractional
differential equation,
%----------------------------------------------------------------
\begin{equation}
    D_{\ast}^{\beta} u(t,x)
    =
    D_{0}^{\alpha} u(t,x), \spc
    t>0, \, 
    x \in \R^N, \, 
    0 < \beta \leq 1, \, 0 < \alpha \leq 2,
    \label{eq:frac_diff_eqn}
\end{equation}
%----------------------------------------------------------------
and then generalize it to distributed order differential equations. 
%----------------------------------------------------------------
 \vspace{0.3cm}

%%%%%%%%%%%%%%%%%%%%%%%%%%%%%%%%%%%%%%%%%%%%%%%%%%%%%%%%%%%%%%%%%
\subsection{Discretization of the time-fractional derivative.}
%%%%%%%%%%%%%%%%%%%%%%%%%%%%%%%%%%%%%%%%%%%%%%%%%%%%%%%%%%%%%%%%%
Using the Caputo time-fractional derivative \cite{Caputo},
the left-hand-side of (\ref{eq:frac_diff_eqn}) becomes
%----------------------------------------------------------------
\begin{equation}
    D_{\ast}^{\beta}u(t,x) =
    \frac{1}{\Gamma(1-\beta)} \int_{0}^t
    \left[ \frac{\partial u(s,x)}{\partial s} \right]
    \frac{ds}{(t-s)^{\beta}}, \,\,
    0 < \beta < 1.
    \label{eq:time_frac_der}
\end{equation}
%----------------------------------------------------------------
Note that when $\beta=1$, 
$D_{\ast}^{\beta}{u(t,x)} = \partial u/ \partial t$.
When $0 < \beta < 1$, we will use the following
discretization (see \cite{LSAT_05} for the derivation):
%----------------------------------------------------------------
\begin{eqnarray}
  D_{\ast}^{\beta} u^n_j
  & \approx &
  \frac{1}{\Gamma(1-\beta)}
  \sum_{m=0}^n \int_{t_n}^{t_{n+1}}
  \frac{ u^{'}_j (t_{n+1}-s) }{ s^{\beta} } ds
  \nonumber \\
  & = & \frac{1}{\nu \tau^{\beta}}
  \left(
     u^{n+1}_j - \sum_{m=1}^n c_m u^{n+1-m}_j - \gamma_n u^0_j
  \right)
  \label{eq:time_frac_der_discr_caputo}
\end{eqnarray}
%----------------------------------------------------------------
where 
%----------------------------------------------------------------
\[
  \gamma_m = (m+1)^{1-\beta} - m^{1-\beta}, \,\, m=0,1,\ldots,n,
  \quad
  c_m = \gamma_{m-1} - \gamma_m, \,\, m=1,\ldots,n
\]
%----------------------------------------------------------------
and $\nu = \Gamma(2-\beta)$.
The formulas for the coefficients 
$c_m$ and $\gamma_m$ and the scalar $\nu$  
that were used in (\ref{eq:time_frac_der_discr_caputo}),  
which were based upon the Caputo time-fractional derivative, easily
generalize to other definitions of the 
time-fractional derivative. For 
example, in the case of the
Grunwald-Letnikov time-fractional derivative,
$\nu = 1$ and $\gamma_m$ and $c_m$ are
re-defined as the following \cite{Cies}:
%----------------------------------------------------------------
\[
  c_m = {\ds \left| 
             \left( \begin{array}{c} \beta \\ m \end{array} \right) 
             \right|, 
        } \, k=1,\ldots,n, 
  \quad
  \gamma_{m} = {\ds 1-\sum_{i=1}^m c_i}, \, m=0,\ldots,n. 
\]
%----------------------------------------------------------------
For simplicity of notation, we will now set 
%----------------------------------------------------------------
\[ \begin{array}{rcl}
  w_0 & = & \gamma_n \\  
  w_i & = & c_{n+i-1}, \,\, i=1,\ldots,n.
\end{array} \]
%----------------------------------------------------------------
and, as a result, (\ref{eq:time_frac_der_discr_caputo}) can be 
rewritten as 
%----------------------------------------------------------------
\begin{equation}
  D_{\ast}^{\beta} u^n_j =
  \frac{1}{\nu \tau^{\beta}}
  \left( u^{n+1}_j - \sum_{m=0}^n w_m u^m_j \right).
  \label{eq:space_time_frac_eqn_w}
\end{equation}
%----------------------------------------------------------------
Note that for $\beta=1$, $\nu=\Gamma(2-\beta)=1$ and 
$w_0=\cdots=w_{n-1}=0$ with $w_n=1$. In this case, 
(\ref{eq:time_frac_der_discr_caputo}) reduces to the
standard forward-time discretization for $\partial u/\partial t$:
%----------------------------------------------------------------
\[
  D_{\ast}^{1} u^n_j =
  \frac{\partial u}{\partial t} \approx
  \frac{u^{n+1}_j - u^n_j}{\tau}.
\]
%----------------------------------------------------------------

%%%%%%%%%%%%%%%%%%%%%%%%%%%%%%%%%%%%%%%%%%%%%%%%%%%%%%%%%%%%%%%%%
\subsection{Discretization of the space-fractional derivative.}
%%%%%%%%%%%%%%%%%%%%%%%%%%%%%%%%%%%%%%%%%%%%%%%%%%%%%%%%%%%%%%%%%
Just as the discretization for the time-fractional derivative assumes a simple
form when $\beta = 1$, the discretization
for the space-fractional 
derivative, based upon centered differences,
assumes a simple form when $\alpha=2$.
For example, when $\alpha=2$ and the $N=2$, 
%----------------------------------------------------------------
\[
  D_{0}^{\alpha} u^n_j =  
  \Delta u^n_j \approx  
  \frac{1}{h^2} \left(
    u^n_{(j_1+1,j_2)} + 
    u^n_{(j_1-1,j_2)} + 
    u^n_{(j_1,j_2+1)} + 
    u^n_{(j_1,j_2-1)} - 4 u^n_{(j_1,j_2)}.
  \right) 
\]
%----------------------------------------------------------------
In $N$-dimensions, the stencil consists of $j=(j_1,\ldots,j_N)$ 
and its nearest $2N$ neighbors with each nearest neighbor 
being $h$ units away from $j$.
When $\alpha=\{\alpha_1,\ldots,\alpha_M\} \neq 2$, the
space-fractional derivative is given by \cite{UmarovSteinberg}:
%----------------------------------------------------------------
\begin{equation}
  D_{0}^{\alpha} u^n_j \approx -q_0(\alpha,h) u^n_j + 
  \sum_{k \neq 0} q_k(\alpha,h) u^n_{j-k}
  \label{eq:space_frac_discr}
\end{equation}
%----------------------------------------------------------------
where the coefficients $q_0(\alpha,h)$ and $q_k(\alpha,h)$ 
are defined in (\ref{eq:q_k}) and (\ref{eq:q_0})
using the multiterm definition for $a(\alpha)$.
The geometric consequence 
of changing $\alpha$ from $\alpha=2$ to 
$\alpha=\{\alpha_1,\ldots,\alpha_M\} \neq 2$ is 
that the stencil gets enlarged from  
$2N+1$ grid points to 
all of the lattice points in $\Z^N_h$.
\vspace{0.3cm}

%%%%%%%%%%%%%%%%%%%%%%%%%%%%%%%%%%%%%%%%%%%%%%%%%%%%%%%%%%%%%%%%% 
\subsection{Construction of the explicit finite difference scheme.} 
%%%%%%%%%%%%%%%%%%%%%%%%%%%%%%%%%%%%%%%%%%%%%%%%%%%%%%%%%%%%%%%%% 
Setting the 
discretizations for the time and space-fractional derivatives equal to each 
other in (\ref{eq:time_frac_der_discr_caputo}) and 
(\ref{eq:space_frac_discr}), we get 
%----------------------------------------------------------------
\begin{equation}
  \frac{1}{\nu \tau^{\beta}}
  \left( u^{n+1}_j - \sum_{m=0}^n w_m u^{m}_j \right) =
  -q_0(\alpha,h) u^n_j +
  \sum_{k \neq 0} q_k(\alpha,h) u^n_{j-k}.
\end{equation}
%----------------------------------------------------------------
Solving for $u^{n+1}_j$, the following
explicit finite-difference scheme is constructed:
%----------------------------------------------------------------
\begin{equation}
  u^{n+1}_{j}
   = 
  {\ds \sum_{m=0}^{n-1} w_m u_j^{m} +
       \sum_{k \in \Z^N} p_k u^n_{j-k}},
  \label{eq:time_space_frac_discr_eqn} 
\end{equation}
%----------------------------------------------------------------
where
%----------------------------------------------------------------
\[
  p_k = \nu \tau^{\beta} Q_k(\alpha,h), \, k \neq 0 
  \quad \mbox{and} \quad
  p_0 = w_n - \nu \tau^{\beta} q_0(\alpha,h).
\] 
%----------------------------------------------------------------
When $\beta = 1$, the coefficients $p_k$ are equivalent to the
transition probabilities $p_k$ in (\ref{eq:markovian_transprob}).
Furthermore, since all the transition probabilities are non-negative
and taking into account that $w_n = c_1 = 2-2^{1-\beta}$ and 
$\nu = Gamma(2-\beta)$, we have an upper bound for the stepsize $\tau$:
%----------------------------------------------------------------
\[
  p_0 \geq 0 
  \quad \Rightarrow \quad 
  0 < \tau \leq 
  \left(  
    \frac{2 - 2^{1-\beta}} {\Gamma(2-\beta) q_0(\alpha,h)}.                  
  \right)^{1/\beta}.
\]
%----------------------------------------------------------------

The update $u^{n+1}_j$ in (\ref{eq:time_space_frac_discr_eqn})
is determined by Markovian contributions
(those values of $u$ at time $t=t_n$) and non-Markovian
contributions (those values of $u$ at times
$t=\{t_0,t_1,\ldots,t_{n-1}\}$).
The order of the time fractional derivative $\beta$
determines the effect that the non-Markovian transition
probabilities ($w_0,\ldots,w_{n-1}$)
has on $u^{n+1}_j$. This effect can be measured by
examining the 
sum of all of the transition probabilities in
(\ref{eq:time_space_frac_discr_eqn}):
%----------------------------------------------------------------
\begin{equation}
  {\ds \sum_{m=0}^{n-1} w_m + 
  \sum_{k \in \Z^N} p_{k} = 1}, 
  \quad
  \left\{ \begin{array}{rcl}
    {\ds \sum_{m=0}^{n-1} w_m} & = & 1 - w_n \\ 
    & & \\
    {\ds \sum_{k \in \Z^N} p_k} & = & w_n.
  \end{array} \right.
  \label{eq:transprob_sum}
\end{equation}
%----------------------------------------------------------------
Recall that when $\beta=1$, $w_n = 1$ and 
$w_0 = \cdots = w_{n-1} = 0$.  In this case,
the first term in (\ref{eq:transprob_sum}) vanishes and  
$p_0 = 1 - \tau q_0(\alpha,h)$.

\begin{figure}[htbp]
  \centering
  \includegraphics[width=5.0in]{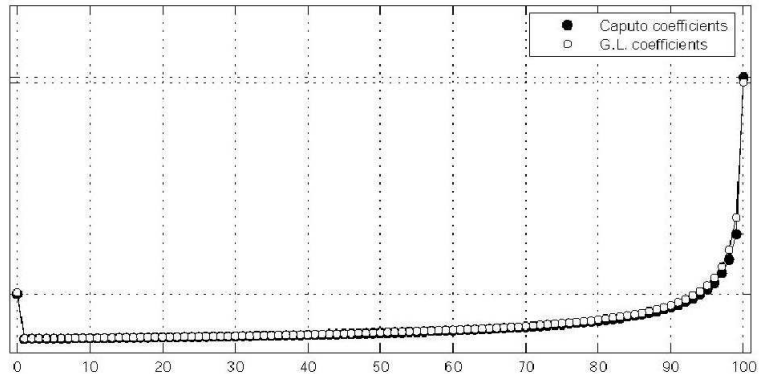}
  \caption{\small The weight $w_m$ associated with the density
   $u_j^m$ is plotted as a function of $m$ for both the Caputo 
   and Gr\"{u}nwald-Letnikov (G.L.) time-fractional derivatives
   and $\beta=0.9$.
   The lower dotted horizontal line corresponds to the value of 
   $w_0 \approx 0.005$ while the upper two dotted lines correspond 
   to $w_n=c_1$ for both the Gr\"{u}nwald-Letnikov 
   ($w_{100}=0.8$) and Caputo derivatives ($w_{100} \approx 0.851$).}
  \label{fig:coef}
\end{figure}

When $0<\beta <1$, the values of $u^n_j$ 
associated with $t \in \{t_0,\ldots,t_{n-1}\}$
are weighted by the coefficients 
$\{w_0,w_1,\ldots,w_{n-1}\}$.
Figure \ref{fig:coef} plots $w_m$ for 
$m=0,1,\ldots,n$ where $n=100$ and $\beta=0.9$.  
It is well-known that the sequence
$\{w_m\}_{m=1}^{n}$ are monotone increasing \cite{Cies}, i.e. 
$w_1 < w_2 < \ldots < w_{n-1} < w_n$.
However, it is not true $w_0 < w_1$.  In fact, 
in Figure \ref{fig:coef}, $w_9 < w_0 < w_8$.  Hence, the 
contribution of $u^0_j$ to $u^{101}_j$ is quite 
large relative to the other intermediate values of $u^n_j$.  
We will see later on that this will have important consequences in 
non-Markovian random walk numerical simulations.   
%----------------------------------------------------------------
%----------------------------------------------------------------

\section{Monte Carlo Protocol for the Random Walk}

%%%%%%%%%%%%%%%%%%%%%%%%%%%%%%%%%%%%%%%%%%%%%%%%%%%%%%%%%%%%%%%%%%
\subsection{General Framework.}
%%%%%%%%%%%%%%%%%%%%%%%%%%%%%%%%%%%%%%%%%%%%%%%%%%%%%%%%%%%%%%%%%%
The random walk model corresponding to the governing equation in
(\ref{eq:frac_diff_eqn}) uses the 
non-Markovian transition probabilities,
$w_m$ the the Markovian transition probabilities $p_k$ to assign
where in the $\Z_h^N$ lattice
a particle will jump to. This jump can be
based upon a partitioning of the unit interval
${\mathcal{P}}=[0,1)$ into two disjoint subintervals
${\mP}_1$ and ${\mP}_2$ such that
${\mP} = {\mathcal{P}}_1 \cup {\mathcal{P}}_2$
where ${\mathcal{P}}_1 = [0,1-w_n)$ and
${\mathcal{P}}_2=[1-w_n,1)$.

We will use a two-dimensional
walk for illustration purposes.
The random walk process begins by
generating a uniformly distributed random number $r$
in the unit interval and observing what subinterval
(${\mathcal{P}}_1$ or ${\mathcal{P}}_2$) it falls into.
If $r \in {\mathcal{P}}_1 = [0,1-w_n)$, then the particle
will do a \emph{non-Markovian} jump, i.e. the jump will be determined
by transition probabilities $w_m, \, m=0,\ldots,n-1$.
Otherwise, if $r \in {\mathcal{P}}=[1-w_n,1)$, then the 
particle will undergo a
\emph{Markovian jump}, i.e. the jump will be determined by transition
probabilities $p_k$.
In effect, 
the random walk interpretation presented here is a two-dimensional extension 
of the one-dimensional random walk interpretation given in \cite{GMMPP}.

\vspace{0.3cm}
\subsection{Non-Markovian Jumps.}
If $0 < \beta < 1$ and 
$r \in {\mathcal{P}}_1$, then the jump 
that the particle takes will be determined by 
$w_m$, $m=0,\ldots,n-1$.
Let 
${\mathcal{A}}=\{{\mathcal{A}}_0,{\mathcal{A}}_1,\ldots,{\mathcal{A}}_{n-1}\}$
be an $n$-element set such that ${\mathcal{A}}_i=w_i$, $i=0,\ldots,n-1$.
Furthermore, let the interval
${\mathcal{P}}_1$ be refined in the following way:
%-----------------------------------------------------------------
\[
  {\mathcal{P}}_1 = [{\mathcal{B}}_0,{\mathcal{B}}_1,\ldots,
  {\mathcal{B}}_n),
\]
%-----------------------------------------------------------------
such that ${\mathcal{B}}_0=0$ and
${\mathcal{B}}_j = \sum_{i=0}^{j-1} {\mathcal{A}}_{i}$, 
$j=1,\ldots,n$.
If $r \in [{\mathcal{B}}_0,{\mathcal{B}}_1)=[0,w_0)$, then
the position of the particle at $t=t_{n+1}$ is given by
$\bd{S}_{n+1}=\bd{S}_0$ (the origin).
Otherwise,
if $r \in [ {\mathcal{B}}_{j-1}, {\mathcal{B}}_{j} )$, $j=1,\ldots,n$,
then the particle will jump back to the position that
it had visited at time $t=t_j$, i.e.
$\bd{S}_{n+1}=\bd{S}_j$.

\vspace{0.3cm}
%%%%%%%%%%%%%%%%%%%%%%%%%%%%%%%%%%%%%%%%%%%%%
\subsection{Markovian Jumps when $\alpha=2$.}  
%%%%%%%%%%%%%%%%%%%%%%%%%%%%%%%%%%%%%%%%%%%%%
If $r \in {\mP}_2=[1-w_n,1)$ and $\alpha=2$ then
the jump will only be to adjacent lattice grid points.
Let ${\mathcal{P}}_2$ be partitioned in the following manner:
%-----------------------------------------------------------------
\[
  {\mP}_1= [{\mB}_0,{\mB}_1,\ldots,{\mB}_5)
\] 
%-----------------------------------------------------------------
where ${\mB}_0 = 1-w_n$ and 
${\mathcal{B}}_j = {\mB}_0 + \sum_{i=0}^{j-1} {\mathcal{A}}_{i}$
$(j=1,\ldots,5)$.
Here, ${\mA} = \{{\mA}_0,{\mA}_1,{\mA}_2,{\mA}_3,{\mA}_4\}$ where 
${\mA}_0 = w_n-4 \eta$ and 
${\mA}_i = \eta = \nu \tau^{\beta}/h^{\alpha}$, $i=1,2,3,4$.  
If $r \in [{\mB}_0,{\mB}_1)$, then the particle remains at the 
current position, otherwise if 
$r \in  
\{ [{\mB}_1,{\mB}_2), 
   [{\mB}_2,{\mB}_3), 
   [{\mB}_3,{\mB}_4), 
   [{\mB}_4,{\mB}_5) \}$
then the particle will move left, right, up or down, 
respectively,
one lattice position.

%%%%%%%%%%%%%%%%%%%%%%%%%%%%%%%%%%%%%%%%%%%%% 
\subsection{Markovian Jumps when 
$\alpha=\{\alpha_1,\ldots,\alpha_M\} \neq 2$.} 
%%%%%%%%%%%%%%%%%%%%%%%%%%%%%%%%%%%%%%%%%%%%% 
If $r \in {\mP}_2=[1-w_n,1)$ 
and $\alpha=\{\alpha_1,\ldots,\alpha_M\} \neq 2$, 
then the jump 
will be determined by an {\em infinite} 
partition refinement of ${\mP}_2$. 
Let \[
  {\mA}=\{{\mA}_0,{\mA}_1,\ldots\}; \spc
  {\mP}_1= [{\mB}_0,{\mB}_1,\ldots)
\]
such that ${\mathcal{B}}_0 = 1-w_n$ and
${\mathcal{B}}_j = {\mB}_0 + \sum_{i=0}^{j-1} {\mathcal{A}}_{i}$ 
$(j=1,\ldots)$.
In this case, the set $\mA$ consists of all of the
transition probabilities $p_k$, $k \in \Z^2$, with 
${\mA}_0=p_0$.
If $r \in [{\mB}_0,{\mB}_1)=[1-w_n,(1-w_n)+p_0)$,
then the particle will remain at the current position. Otherwise,
if $r \in [{\mB}_s,{\mB}_{s+1})$,
then there exists a unique $k=(k_1,k_2) \in \Z^2$
associated with $s \in \N$
such that
the particle will jump from
$\bd{S}_n$ to $\bd{S}_{n+1} = \bd{S}_n + (k_1 h, k_2 h)$.

\section{Simulations}

\begin{figure}[htbp]
  \centering
  \includegraphics[width=3.0in]{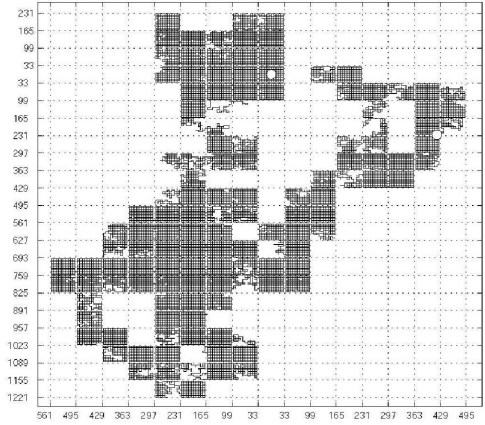}
  \caption{\small This random walk simulation depicts
  classical Brownian motion confined to rectangular 
  cytosketetal barriers.  The parameters used in this simulation 
  are as follows:
  $h=6$ nanometers, $\tau=10^{-6}s$ and 
  $a=9 \times 10^{-12} \mbox{m}^2/s$.  The 
  barriers are spaced out every 66 nanometers and the the probability 
  of escape is $p=0.01$ when a particle encounters a barrier.}
  \label{fig:kusumi}
\end{figure}

Our motivation of the numerical simulations presented here is to see how DODE 
simulations of biomolecular motion of particles on a cell surface differ from 
those based upon classical Brownian motion. Although the DODE random walk 
models are described theoretically for multivariate case in $N$-dimensions, 
nevertheless all our simulations are conducted in the two dimensional case 
since we are interested in the diffusion of proteins on a cell membrane 
surface, which can be locally approximated by a two-dimensional membrane 
sheet. In \cite{Kusumi2}, simulated particle motion is based upon the 
classical Brownian motion scenario (where $\alpha = 2$ and $\beta = 1$) in 
which the particle is confined within cytoskeletal barriers (see Figure 
\ref{fig:kusumi}). In these single particle tracking studies, particle 
appears to be spatially and temporarily confined within \emph{transient 
confinement zones}.  Although the barriers are never directly observed, it is 
postulated that the cytoskeletal barriers are the reason for the transient 
spatial confinement of particle. In principle, DODE simulations provide an 
alternative explanation for the observed trajectories in single particle 
tracking studies that does not necessarily require the existence of 
cytoskeletal barriers to explain transient confinement.

In \cite{Kusumi2}, the authors use the mean-squared-displacement 
formula $4 a \tau = h^2$ in which the parameters $a$ (the diffusion 
coefficient) , $\tau$ (the timestep) and $h$ (the lattice width), 
respectively, are given using the following values:
$h=6$ nanometers and $\tau = 1 \mu s$ (microseconds, or
$\tau=10^{-6}$ seconds).  Since the 
mean-squared displacement formula implicitly assumes that
%------------------------------------------------------------------------
\[
  p_0 = 1 - 4 a \frac{\tau^{\beta}}{h^{\alpha}} = 
        1 - 4a \frac{\tau}{h} = 0,
\]
%------------------------------------------------------------------------
the diffusion coefficient is then computed as 
$a = h^2/(4\tau) = 9 \times 10^{-12} \mbox{m}^2/s$. 
To facilitate a comparison of our DODE simulations with the simulations
of \cite{Kusumi2,Kusumi}, we will also use the same diffusion coefficient
($a_1 = \cdots = a_M = a = 9 \times 10^{-12} 
\mbox{m}^2/s$) and the same lattice width
($h=6$ nanometers).
Using the fact that the transition probabilities sum to 1,
%------------------------------------------------------------------------
\[
  1 = \sum_{m=0}^{n-1} w_m + \sum_{k} p_k =
  (1-w_n) + p_0 +
  \nu \tau^{\beta} q_0(\alpha,h)
\]
%------------------------------------------------------------------------
we can now solve for $\tau$ in terms of
$\alpha$, $\beta$ and $p_0$,
%------------------------------------------------------------------------
\[
  \tau = \tau(\alpha,\beta,p_0) =
  \left( \frac{c_1 - p_0}{\nu q_0(\alpha,h)}
         \right)^{1/\beta} = 
  \left(
    \frac{(2-2^{1-\beta}) - p_0}{\Gamma(2-\beta) q_0(\alpha,h)}
  \right)^{1/\beta}. 
\]
%------------------------------------------------------------------------
As in \cite{Kusumi2}, we set
$p_0 = 0$.
However, due to the dependence of $\tau$ on $\alpha$ and $\beta$,
the relative
size of the timestep (from $\tau=10^{-6}s$ in the case of $\alpha=2$ and 
$\beta=1$)
will change as $\alpha$  and $\beta$ vary.
Instead if fixing the simulations to have the same stepsize $\tau$, we will fix
the duration of the overall walk to be the same,
Let $T$ denote the overall duration of the random walk
simulation.  In all of our DODE simulations, $T$ is set to 
$T = \frac{1}{30}$ seconds.  This is equivalent to 1 frame at video 
rate where video rate is
measured as 30 frames per second.  All simulations were performed in 
MATLAB\cite{matlab}.

%%%%%%%%%%%%%%%%%%%%%%%%%%%%%%%%%%%%%%%%%%%%%%%%%%%%%%%%%%%%%%%%%%%%%%%%
%%%%%%%%%%%%%%%%%%%%%%%%%%%%%%%%%%%%%%%%%%%%%%%%%%%%%%%%%%%%%%%%%%%%%%%%

%%%%%%%%%%%%%%%%%%%%%%%%%%%%%%%%%%%%%%%%%%%%%%%%%%%%%%%%%%%%%%%%%%%%%%%%
\begin{figure}[ht]
  \centering  
  \includegraphics[width=5.5in]{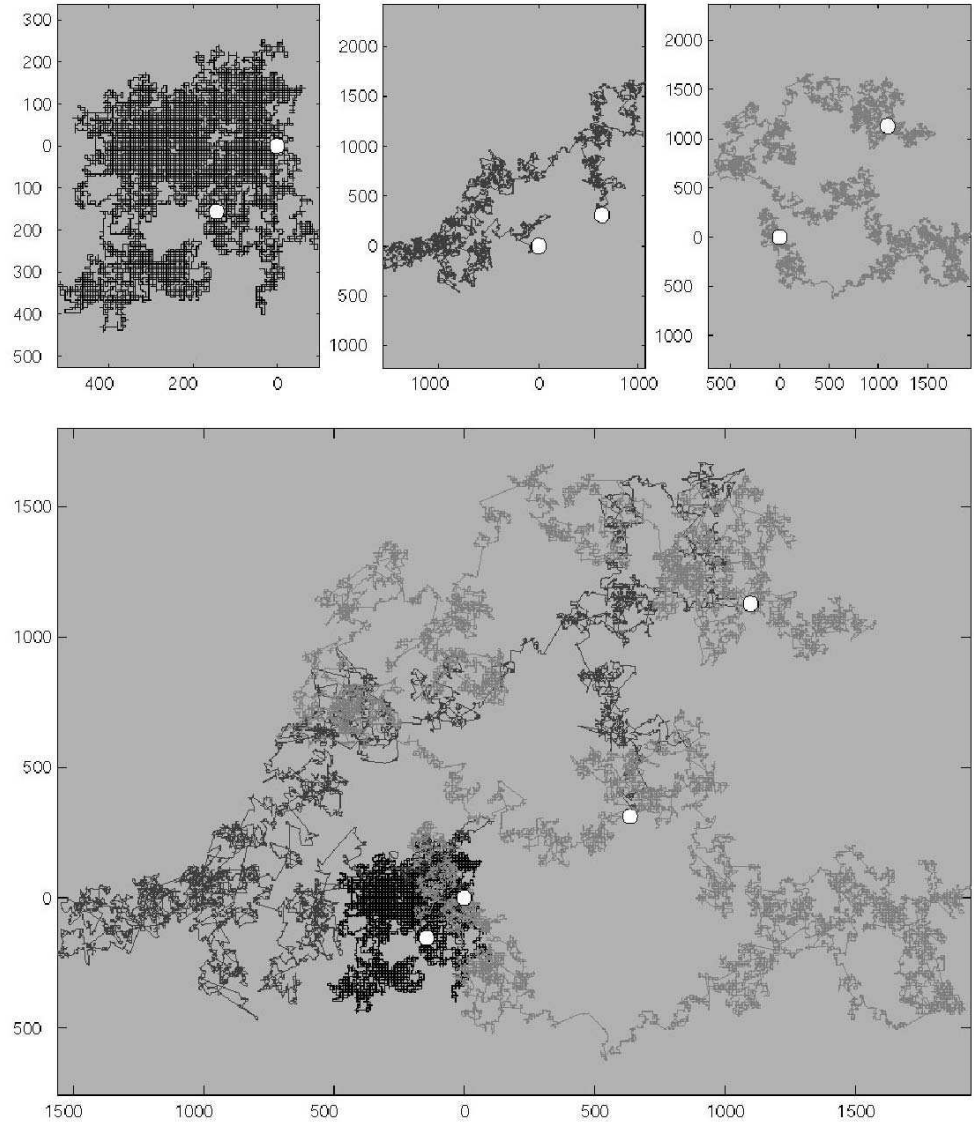}
  \caption{\fns The first three subplots in the top row correspond to 
           Markovian DODE simulations ($\beta=1$) with different values of
           $\alpha$:  $\alpha=2$, $\alpha=1.5$ and $\alpha=\{1.5,2\}$ for 
           the left, middle and right plots.  The bottom plot 
           superimposes all of the top three simulations on one graph.} 
  \label{fig:plot1}
\end{figure}
%%%%%%%%%%%%%%%%%%%%%%%%%%%%%%%%%%%%%%%%%%%%%%%%%%%%%%%%%%%%%%%%%%%%%%%%

%%%%%%%%%%%%%%%%%%%%%%%%%%%%%%%%%%%%%%%%%%%%%%%%%%%%%%%%%%%%%%%%%%%%%%%%
\begin{figure}[ht]
  \centering
  \includegraphics[width=5.5in]{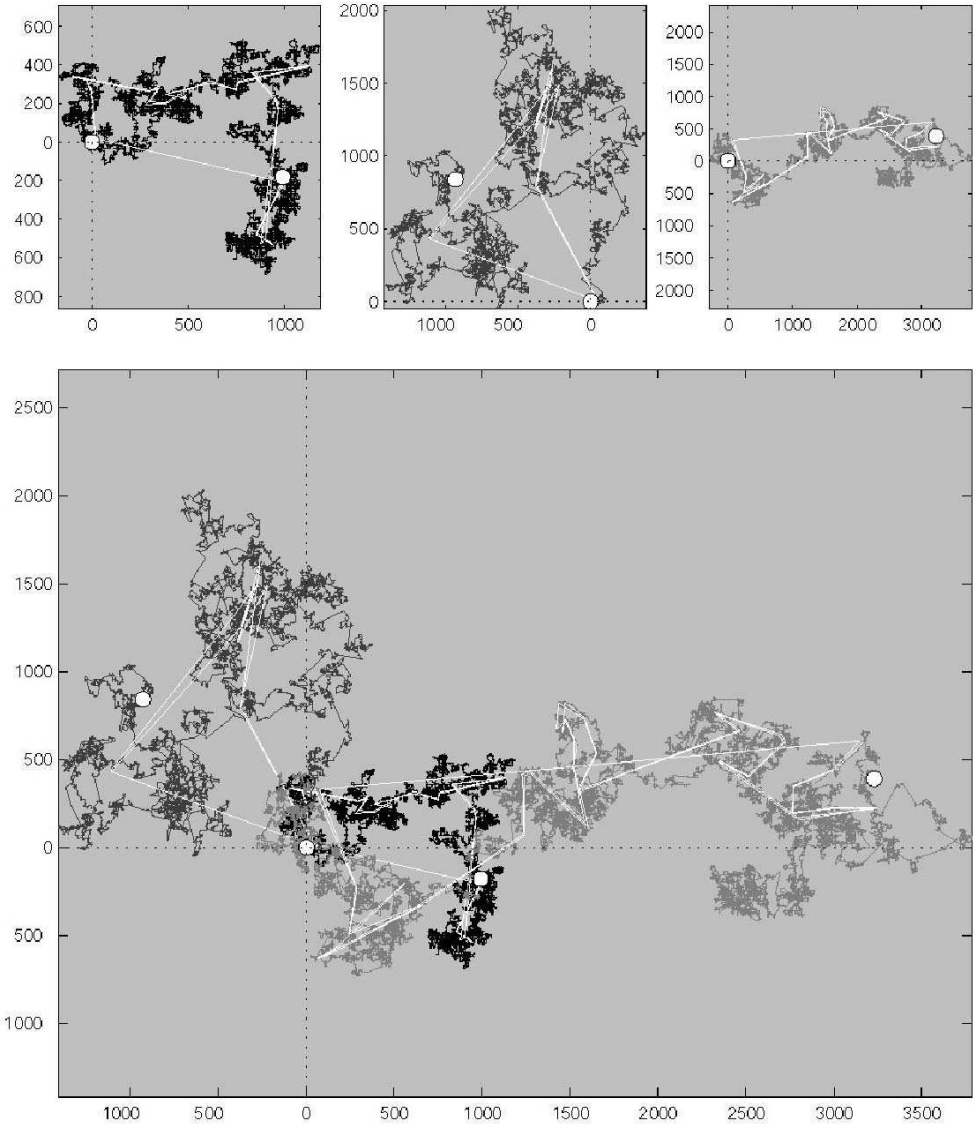}
  \caption{\fns The first three subplots in the top row correspond to
           non-Markovian DODE simulations ($\beta=0.999$) with different 
           values 
           of $\alpha$:  $\alpha=2$, $\alpha=1.5$ and $\alpha=\{1.5,2\}$ 
           for
           the left, middle and right plots.  The dark shaded lines 
           correspond to non-Markovian walks while the white lines
           indicate non-Markovian jumps to previously visited positions.
           The bottom plot
           superimposes all of the top three simulations on one graph.
          }
  \label{fig:plot2}
\end{figure}
%%%%%%%%%%%%%%%%%%%%%%%%%%%%%%%%%%%%%%%%%%%%%%%%%%%%%%%%%%%%%%%%%%%%%%%%

%%%%%%%%%%%%%%%%%%%%%%%%%%%%%%%%%%%%%%%%%%%%%%%%%%%%%%%%%%%%%%%%%%%%%%%%
\begin{figure}[hb]
  \centering
  \includegraphics[width=5.5in]{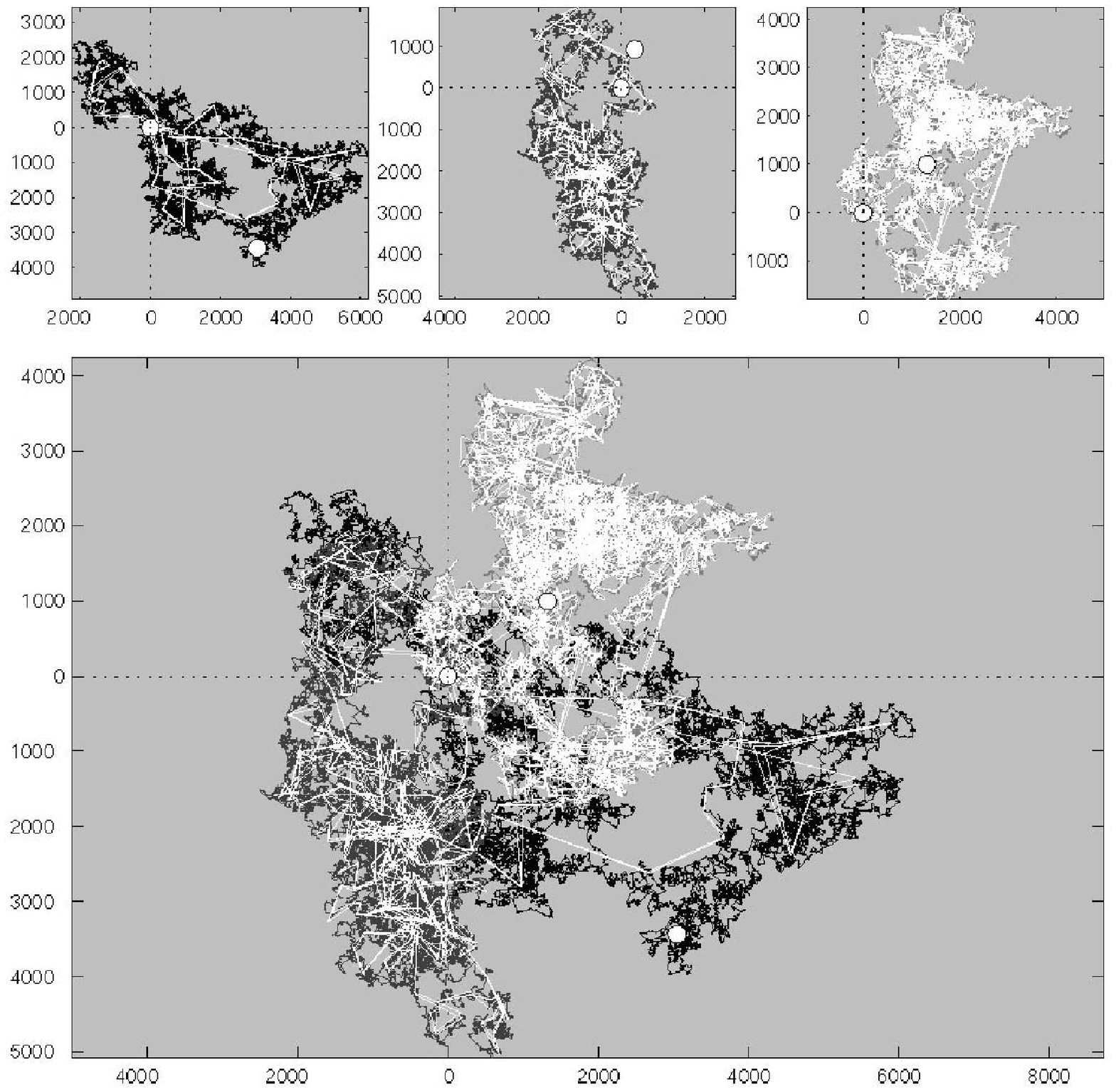} 
  \caption{\fns The first three subplots in the top row correspond to
           non-Markovian DODE simulations with $\alpha=\{0.8,1.3,1.8\}$
           and different values of $\beta$:  
           $\beta=.999$, $\beta=.99$ and $\beta=.999$
           for
           the left, middle and right plots.  The dark shaded lines
           correspond to non-Markovian walks while the white lines
           indicate non-Markovian jumps to previously visited positions.
           The bottom plot
           superimposes all of the top three simulations on one graph.
  }
  \label{fig:plot3}
\end{figure}
%%%%%%%%%%%%%%%%%%%%%%%%%%%%%%%%%%%%%%%%%%%%%%%%%%%%%%%%%%%%%%%%%%%%%%%%

Figure \ref{fig:plot1} shows various Markovian DODE 
simulations ($\beta=1$) across 
various values of $\alpha$.  The left, middle and right plots in the 
top row show DODE simulations
for $\alpha = \{ 2 \}$, $\alpha=\{ 1.5 \}$ and
$\alpha = \{ 1.5, 2 \}$, respectively. 
The first two DODE simulations
are actually monofractal DODE simulations with $M=1$ while the last one 
($\alpha=\{1.5,2\}$) is a multi-fractal case with $M=2$.
The large white dots indicate the first and last positions 
of the random walk and the starting position is 
always the origin $(0,0)$. 
It is clear that for these DODE simulations with
$\alpha \not= \{2\}$ that the particle
travels much longer distances since
the probability of jumping to faraway lattice sites is greater
than what would be expected for $\alpha=2$.

Figure \ref{fig:plot2} shows  
various non-Markovian DODE simulations ($\beta=0.999$) using the same values of 
$\alpha$ as in Figure \ref{fig:plot1}.  The bottom plot in both Figures
\ref{fig:plot1}  and \ref{fig:plot2} show the plots on top row superimposed 
on one graph.  The dark shaded lines correspond to Markovian jumps 
($r \in {\mP}_1$)
while the 
white lines correspond to non-Markovian
jumps ($r \in {\mP}_2$).
The frequency of the non-Markovian jumps are given by the size of the 
${\mP}_1$ interval.  For $\beta=0.999$, 
${\mP}_1 \approx [0,1-w_n) = [0,0.00069339)$.  
Hence, the probability at every 
timestep of doing a non-Markovian jump is $0.00069339$.
The bottom plot in Figure \ref{fig:plot2} shows the superposition all three
non-Markovian DODE simulations on the same graph. 

For Figure \ref{fig:plot3}, we have non-Markovian DODE simulations for 
a fixed set of $\alpha$ values ($\alpha=\{0.8,1.3,1.8\}$) with 
$\beta$ varying.  The left, middle and right plots correspond to 
$\beta=0.999$, $\beta=0.99$  and $\beta=0.9$, respectively.  
The probability of taking a non-Markovian per timestep for these graphs 
is $0.00069339$ (left), $0.0070$ (middle) and $0.0718$ (right).
For example, roughly 7\% of all jumps for the right subplot on the top row
are non-Markovian jumps. The effect of decreasing $\beta$ is clear:  the 
overall 
distances that the particle traverses is decreased since motion is
constrained by jumps to previously visited positions.

The average jump sizes  
associated with Figures \ref{fig:plot1}, \ref{fig:plot2} 
and \ref{fig:plot3} are shown in 
Table \ref{tab:avg_jump}.  The numbers in the brackets before the colon 
correspond to the $(\alpha,\beta)$ pair used in the DODE simulation 
while the number after the colon corresponds to the average jump size.
For the non-Markovian walks, the average jump 
length is larger when, for a fixed set of $\alpha$ values, $\beta$ is decreased 
from 1. This is a consequence of the non-Markovian nature of the random walks for 
$0 < \beta < 1$.  
Since the particle is allowed to jump back to any previously visited 
position, the jump size can be quite large if the previously visited position was 
spatially remote from the particle's 
current position (see Figure \ref{fig:plot3}).
In particular, in Figure \ref{fig:coef}, the probability of the particle to 
jump back to the origin is disproportionately larger than for other previously 
visited sites.  In Figures \ref{fig:plot2} and \ref{fig:plot3}, one can 
observe evidence of this phenomenon.    

%%%%%%%%%%%%%%%%%%%%%%%%%%%%%%%%%%%%%%%%%%%%%%%%%%%%%%%%%%%%%%%%%%%%%%%%
\begin{table}[h]
\centering
\caption{This table reports the average jump size (after the colon)
for all of the DODE 
simulations in Figures \ref{fig:plot1}, \ref{fig:plot2} and
\ref{fig:plot3}.  The numbers before the colon indicates values of the
$(\alpha,\beta)$-pair used in the DODE simulation.}
\begin{tabular}{|l|l|l|l|} \hline
& 
{\fns Left Plot} & 
{\fns Middle Plot} & 
{\fns Right Plot} \\ \hline
{\fns Figure \ref{fig:plot1}} & 
       {\tiny (2,1): 6.0000} & 
       {\tiny (1.5,1): 10.9770} & 
       {\tiny (\{.5,2\},1): 7.3320} \\ \hline
{\fns Figure \ref{fig:plot2}} & 
         {\tiny (2,0.999): 6.0038} & 
         {\tiny (1.5,0.999): 11.0707} & 
         {\tiny (\{1.5,2\},0.999): 7.3593} \\ \hline
{\fns Figure \ref{fig:plot3}}  & 
         {\tiny (\{0.8,1.3,1.8\},0.999): 17.0328} & 
         {\tiny (\{0.8,1.3,1.8\},0.99):  17.1663} & 
         {\tiny (\{0.8,1.3,1.8\},0.9):   19.8946} \\ \hline 
\end{tabular}
\label{tab:avg_jump}
\end{table}

\section{Conclusion}

Qualitatively, the DODE simulations provide a richer repertoire of motion, 
compared to monofractal walks when $M=1$. Macroscopically, the DODE 
trajectories tend to cluster together more often than the monofractal 
walks.  The clustering is even more pronounced when the motion is 
non-Markovian due to the memory the particle has for previously visited 
positions. Moreover, one does not have to hypothesize the existence of 
barriers to explain why a particle appears trapped in a transient 
confinement zone or hops large distances. The clustering of trajectories 
and large jumps are a natural consequence of the DODE random walk model. 
However, when the motion is non-Markovian, the particle has a strong 
propensity to jump back to the origin, a consequence of the 
disproportionately large weight $w_0$ associated with $u^0_j$. While 
jumping back to previously visited ``compartments'' is observed for 
experimentally observed single particle tracking data \cite{Kusumi}, one 
does not experimentally observe molecules jumping back from its current 
position to the starting point.  Nonetheless, the DODE random walk models 
closely resemble the data from single particle tracking experiments of 
molecules moving on cell membranes\cite{Kusumi2,Kusumi}.  This is not 
surprising since the motion of biomolecules on the cell surface occurs in a 
very heterogeneous environment.

\vskip 1cm
%%%%%%%%%%%%%%%%%%%%%%%%%%%%%%%%%%%%%%%%%%%%%%%%%%%%%%%%%%%%%%%%%%%%%%%%%%%%%%%%
\end{document}